\documentclass[12pt,times new roman]{article}
\usepackage{amsfonts}
\usepackage{amssymb}
\usepackage{amsmath}
\usepackage{amscd}
\usepackage{tikz-cd}
\usepackage{times}
\usepackage{amsthm}
\usepackage{hyperref}
\usepackage[margin=1.5in]{geometry}
\usepackage{multirow}
\usepackage{indentfirst}
\usepackage{dsfont}
\usepackage{setspace}
\newtheorem{Theorem}{Theorem}[section]
\newtheorem{Lemma}[Theorem]{Lemma}

\newtheorem{remark}[Theorem]{Remark}
\newtheorem{innercustomthm}{Theorem}

\newcommand{\floor}[1]{\lfloor #1 \rfloor}
\newenvironment{mythm}[1]
  {\renewcommand\theinnercustomthm{#1}\innercustomthm}
  {\endinnercustomthm}

\begin{document}
\date{September 8th, 2025}
\author{Zeyu Shen}

\title{Higher $G$-theory of simplicial toric varieties and 
vanishing of Chow groups}
\maketitle

\begin{abstract}
This paper gives computations of all the
$G$-theory groups of several classes of simplicial toric 
varieties, including all affine toric surfaces 
when the base field is algebraically closed and has 
characteristic zero, all weighted projective spaces over any 
field and  resolution of singularities of all affine toric 
surfaces $\operatorname{Spec}(k[x,xy,xy^2,...,xy^d])$ over any 
field $k$. The $G$-theory groups 
$G_0,G_1,G_2$ are computed for the product of any two weighted 
projective spaces over any field.
The dimension of the rational vector space 
$G_0(X)\otimes\mathbb{Q}$ for any complete, simplicial toric 
variety $X$ over an algebraically closed field of 
characteristic zero
is shown to be equal to the sum of the Betti numbers of even 
degrees. We also prove that the Chow group $A^2(X)$ of 
codimension 2 cycles
vanishes for any affine, smooth toric variety,
thereby proving a special case of the conjecture that the order 
of this Chow group $A^2(X)$ divides the determinant of the 
matrix whose columns are the minimal generators of the fan 
for any affine, simplicial toric variety $X$.
\end{abstract}
\section{Introduction}
Quillen's algebraic $G$-theory of Noetherian schemes 
provides a natural extension of algebraic $K$-theory of 
Noetherian, separated schemes. It is particularly suited for 
the study of singular schemes. Toric varieties form an 
important testing ground for the algebraic $G$- and $K$-theory 
of Noetherian schemes: their combinatorial structure allows 
explicit computations, and at the same time they capture deep 
phenomena in intersection theory and motivic cohomology.\par
The algebraic $K$-theory of toric varieties has been studied 
extensively by Morelli, Gubeladze, Cortiñas et 
al. in \cite{M},\cite{G2},\cite{CHWW}.
Other authors such as Joshua and Krishna, have considered the 
equivariant algebraic $G$-theory of toric stacks in \cite{JK}.
However, the higher algebraic $G$-theory of affine and 
projective simplicial toric varieties and $G_0$ of complete, 
simplicial toric varieties remain much less understood. 
In particular, there appears to be no complete computations of 
$G$-theory groups for basic classes of simplicial toric 
varieties, 
nor structural results describing $G_0$ in the complete case.
The purpose of this paper is to fill this gap by
explicit computations and new structural results for the $G$-theory of simplicial toric varieties.
The results in this paper provide both explicit data and 
structural insights into the $G$-theory of toric varieties, 
extending previous work on their $K$-theory and supplying 
evidence for conjectural connections between $G$-theory and 
algebraic cycles. \par
In Section 2, we compute all $G$-theory groups of affine 
toric surfaces and resolution of singularities for affine 
toric surfaces $\operatorname{Spec}(k[x,xy,xy^2,...,xy^d])$. We 
also verify a special case of the 
conjecture in \cite{S1}. \par
Recall that for any field $k$, $G_0(k)\cong\mathbb{Z}$.\par
\begin{mythm}{2.6}
Let $X$ be an affine toric surface over an algebraically closed 
field $k$ of characteristic zero.
Let $\delta$ be the determinant of the matrix whose columns are 
the minimal generators of the fan of $X$.
Then $G_n(X)\cong G_n(k)$ for all odd $n$ and $G_n(X)\cong G_n(k)\oplus\mathbb{Z}/\delta\mathbb{Z}$ for all even $n$.
\end{mythm}
\begin{mythm}{2.8}
Let $k$ be any field and let $X=\operatorname{Spec}(k[\sigma^{\vee}\cap\mathbb{Z}^2])$, where $\sigma$ is the cone 
in $\mathbb{R}^2$ generated by $e_2,de_1-e_2$ and $d$ is any 
positive integer. Let $\tilde{X}$ be the resolution of 
singularities of $X$ obtained by adding the ray generated by 
$e_1$ to the fan of $X$. Then $G_n(\tilde{X})\cong G_n(k)^{\oplus 2}$ for all $n$.
\end{mythm}
\begin{mythm}{2.9}
Let $X$ be any affine, smooth toric variety. The Chow group of 
codimension 2 cycles of $X$ vanishes: $A^2(X)=0$.
\end{mythm}
If we combine Theorem 2.6 and Theorem 2.9, we see that 
$G_0(X)\cong\mathbb{Z}$ for any affine smooth toric variety 
$X$, which recovers a special case of Theorem 3.3 in 
\cite{S1}.\par
In Section 3, we compute the dimension of the rational 
vector space $G_0(X)\otimes\mathbb{Q}$, where $X$ is any 
complete, simplicial toric variety over an algebraically closed 
field of characteristic zero. \par
\begin{mythm}{3.1}
Let $X$ be any complete, simplicial toric variety over an 
algebraically closed field of characteristic zero. The 
dimension of the rational vector space 
$G_0(X)\otimes\mathbb{Q}$ is the sum of the Betti numbers of $X$
of even degrees.
\end{mythm}
In Section 4, we compute all the $G$-theory groups of any 
weighted projective 
space over an arbitrary field, improving Joshua and Krishna's 
result in \cite{JK}, which only works for algebraic $G$-theory 
of weighted projective spaces using rational coefficients.
\begin{mythm}{4.2}
Let $k$ be any field and let $X$ be any $d$-dimensional 
weighted projective space over the field $k$. Then $G_n(X)\cong G_n(k)^{\oplus (d+1)}$ for all $n$.
\end{mythm}
In section 5, we compute the first few $G$-theory groups 
$G_0,G_1,G_2$ of a product of two weighted projective spaces 
over any field.\par
\begin{mythm}{5.1}
Let $X, Y$ be any two weighted projective spaces over a field 
$k$. Then for $n=0,1,2$, we have
$G_n(X\times_k Y)\cong\bigoplus_{i+j=n}[G_{i}(X)\otimes_{\mathbb{Z}}G_j(Y)]$.\par
\end{mythm}

\section{Higher $G$-theory of affine toric surfaces and related 
results}
In order to compute all the $G$-theory groups of affine toric 
surfaces over an algebraically closed field of characteristic 
zero, we need to establish several lemmas first. \par
\begin{Lemma}
Let $k$ be an algebraically closed field of characteristic 
zero and $R=k[x,xy,xy^2,...,xy^m]$ 
be the subring of the polynomial ring $k[x,y]$, 
where $m$ is a positive integer.
The image of $x\in G_1(R[\frac{1}{x}])$ under the boundary map $\partial: 
G_1(R[\frac{1}{x}]) \rightarrow G_0(R/xR)$ is $m$.
\end{Lemma}
\begin{proof}
Let $\sigma$ be the cone in $\mathbb{R}^2$ generated by $e_1,e_1+me_2$.
Then $R$ is the semigroup algebra 
$k[\sigma\cap\mathbb{Z}^2]$.
Let $v$ denote $xy^m$.
Note that the elements of the semigroup $\sigma\cap\mathbb{Z}^2$ form a $k$-basis for the $k$-algebra 
$R$.
Since multiplication by $x$ in the ring $R$ corresponds to 
right 
translation by one unit for points in $\sigma\cap\mathbb{Z}^2$, 
there is a  
$k$-basis for $R/xR$ corresponding to the points in $\sigma\cap\mathbb{Z}^2$
which lie outside the cone $\sigma$ when translated one unit to the left.
The elements of this $k$-basis for $R/xR$ are of the form $v^l.1,v^l.
(xy),v^l(xy^2),...,v^l(xy^{m-1})$, where $l=0,1,2,...$.
Hence, as a $k[v]$-module, we have 
\[ \displaystyle R/xR\cong k[v].1\oplus k[v](xy)\oplus k[v]
(xy^2)\oplus...\oplus k[v](xy^{m-1}) (*)\]
Note that for $p=1,2,...,m-1$, we have $(xy^p)^m=(x^{m-p})(xy^m)^p=0\in 
R/xR$.
Hence, the nilradical of $R/xR$ is generated by the quotient images of 
$xy,xy^2,...,xy^{m-1}$ in $R/xR$.\\
Therefore, the ring $(R/xR)_{red}\cong k[xy^m]=k[v]\cong k[t]$.\\
Hence, $(*)$ above implies that the image of $x$ in $\mathbb{Z}\cong 
G_0((R/xR)_{red})$ is $m$, since the ring $(R/xR)_{red}\cong k[t]$ is a 
principal ideal domain so that the class $[R/xR]=[k[v].1\oplus k[v].
(xy)\oplus\ k[v](xy^2)\oplus ...\oplus k[v].(xy^{m-1})]\in 
G_0((R/xR)_{red})$ corresponds to $m$.\\
By dévissage, we know that $G_0(R/xR)\cong G_0((R/xR)_{red})$.\\
Therefore, we conclude that the image $\partial(x)=[R/xR]\in 
G_0(R/xR)$ is $m$.
\end{proof}
\begin{Lemma}
Let $k$ be an algebraically closed field of characteristic zero.
Let $m>0$ be an integer and let $X=\operatorname{Spec}
(k[x,xy,xy^2,...,xy^m])$.
Then $G_n(X)\cong G_n(k)\oplus \mathbb{Z}/m\mathbb{Z}$ for 
every non-negative even number $n$, and that $G_n(X)\cong 
G_n(k)$ for every positive odd number $n$.
\end{Lemma}
\begin{proof}
Since $G_1(k)\cong G_1((R/xR)_{red})\cong G_1(R/xR)$, $G_1(k)$ is mapped 
into $G_0(R/xR)$ by the composition of maps $G_1(R/xR)\rightarrow 
G_1(R)\rightarrow G_1(R[\frac{1}{x}])\rightarrow G_0(R/xR)$ in the $G$-theory 
localization sequence of $R$ and its localization $R[\frac{1}{x}]$.\\
The exactness of the $G$-theory localization sequence implies that this 
composition of maps is zero.\\
Hence, the image of $G_1(k)$ under the boundary map $\partial$ is zero.\\
Similar reasoning implies that the image of $G_n(k)$ under the boundary map
$\partial:G_n(R[\frac{1}{x}])\rightarrow G_{n-1}(R/xR)$ is zero for every positive 
integer $n$.                                  (1)\\
Therefore, Lemma 2.1 above implies that the image of 
$\partial:G_1(R[\frac{1}{x}])\rightarrow G_0(k)$ is $m\mathbb{Z}$.\\
By exactness of the $G$-theory localization sequence for $R$
and its localization $R_x$ we have $G_0(R)\cong \mathbb{Z}\oplus 
\mathbb{Z}/m\mathbb{Z}$.\par
Now we compute the higher $G$-theory of $R$.
For $n=1$, the relevant part of the $G$-theory localization sequence 
for $R$ and its localization $R[\frac{1}{x}]$ is 
\[G_2(R_x)\rightarrow 
G_1(R/xR)\rightarrow G_1(R)\rightarrow G_1(R[\frac{1}{x}])\rightarrow 
G_0(R/xR)\rightarrow ...\]\\
For every $\alpha\in G_1(k)$, we have $\partial(\alpha\cup 
x)=\alpha\partial(x)=m\alpha$.\\
Since $k$ is an algebraically closed field of characteristic zero, it 
contains $m$ distinct $m$-th roots of unity.
Hence the map $G_1(k)=k^*\xrightarrow {m} G_1(k)$ is surjective.\\
(1) above implies that the boundary map sends the $G_2(k)$ direct summand 
of $G_2(R[\frac{1}{x}])\cong G_2(k)\oplus G_1(k)$ to zero.
And this boundary map restricts to the multiplication by $m$ map on the 
$G_1(k)$ direct summand of $G_2(R[\frac{1}{x}])$.
Therefore, the image of this boundary map is $G_1(k)$.\\
The exactness of the $G$-theory localization sequence of $R$ 
and 
its localization $R[\frac{1}{x}]$ implies that the map $G_1(R/xR)\rightarrow G_1(R)$ \
is zero, so that the map $G_1(R)\rightarrow G_1(R[\frac{1}{x}])$ is injective.\\
Hence, $G_1(R)$ is isomorphic to the image of the map 
$G_1(R)\rightarrow G_1(R[\frac{1}{x}])$.\\
By the exactness of the $G$-theory localization sequence for $R$ and $R[\frac{1}{x}]$, 
this image equals the kernel of the map $\partial:G_1(R[\frac{1}{x}])\rightarrow 
G_0(k)$.\\
(1) above implies that this kernel is equal to $G_1(k)$.
Therefore, we deduce that $G_1(R)\cong G_1(k)$.\par
For the case $n=2$, consider the relevant part of the $G$-theory 
localization sequence of $R$ and its localization $R[\frac{1}{x}]$, which is 
\[G_3(R[\frac{1}{x}])\rightarrow G_2(k)\rightarrow G_2(R)\rightarrow 
G_2(R[\frac{1}{x}])\rightarrow G_1(k)\rightarrow...\]\\
Similar to the case $n=1$ above, $\partial(\alpha\cup x)=m
\alpha$ for every $\alpha\in G_2(k)$.
Since $k$ is an algebraically closed field of characteristic zero, 
$G_2(k)$ is a uniquely divisible abelian group.
So $G_2(k)\xrightarrow {m} G_2(k)$ is surjective.
Hence,by a reasoning similar to the $n=1$ case above, the image of 
$\partial:G_3(R_x)\rightarrow G_2(k)$ is $G_2(k)$.
By the exactness of the $G$-theory localization sequence for $R$ and its 
localization $R_x$, $G_2(R)\cong ker(\partial:G_2(R_x)\rightarrow 
G_1(k))$, which is equal to $G_2(k)\oplus ker(G_1(k)\xrightarrow {m} G_
1(k))$.\\
Hence, $G_2(R)\cong G_2(k)\oplus \mathbb{Z}/m\mathbb{Z}$.\par
By Theorem 1.6 in Chapter VI of the book \cite{W} above, 
$G_n(k)$ is a uniquely 
divisible abelian group for positive even numbers $n$, and 
$G_n(k)$ is the direct sum of a uniquely divisible abelian 
group and $\mathbb{Q}/\mathbb{Z}$ for all positive odd numbers 
$n$.
Since the group $\mathbb{Q}/\mathbb{Z}$ is isomorphic to the 
group of 
roots of unity in $\mathbb{C}$, the multiplication by $m$ endomorphism 
of $\mathbb{Q}/\mathbb{Z}$ is surjective. And its kernel is the group 
of $m$-th roots of unity of $\mathbb{C}$, which is cyclic of order 
$m$.
Therefore, by a reasoning similar to the cases $n=1,n=2$ above, we get 
that $G_n(R)\cong G_n(k)\oplus \mathbb{Z}/m\mathbb{Z}$ for every 
positive even number $n$, and that $G_n(R)\cong G_n(k)$ for every 
positive odd number $n$.
Recall that the $G$-theory of 
$X=\operatorname{Spec}(R)$ coincides with the $G$-theory of 
$R$. Therefore, we deduce that $G_n(X)\cong 
G_n(k)\oplus \mathbb{Z}/m\mathbb{Z}$ for every non-negative 
even number $n$, and that $G_n(X)\cong G_n(k)$ for every 
positive odd number $n$.
\end{proof}
\begin{Lemma}
Suppose $\sigma$ is the cone in $\mathbb{R}^2$ generated 
by $e_1,ae_1+e_2$ 
for some integer $a$. Let $k$ be an algebraically closed 
field of characteristic zero, and let $X$ be 
$\operatorname{Spec}(k[\sigma\cap\mathbb{Z}^2])$. Then 
$G_n(X)\cong G_n(k)$ for every $n\geq 0$.
\end{Lemma}
\begin{proof}
Note that left multiplication by the matrix 
$\begin{pmatrix}
    1 & -a\\
    0 & 1
\end{pmatrix}$\
maps $(1,0)$ to $(1,0)$ and $(a,1)$ to $(0,1)$.\\
Since the matrix is in $GL(2,\mathbb{Z})$, left 
multiplication by this matrix represents an automorphism 
of the lattice $\mathbb{Z}^2$.
Hence, we deduce that the semigroup algebra 
$k[\sigma\cap\mathbb{Z}^2]$
is isomorphic to the ring $k[x,y]$.\\
By the fundamental theorem of $G$-theory of Noetherian rings, 
we have $G_n(X)\cong G_n(k[x,y])\cong G_n(k)$ for 
every integer $n\geq 0$.
\end{proof}
\begin{remark}
Let $a,b>0$ be two relatively prime integers and $\sigma$ 
be the cone in 
$\mathbb{R}^2$ generated by $e_1,ae_1+be_2$.
We have computed the $G$-theory of $k[\sigma\cap\mathbb{Z}^2]$ 
for the case $b=1$.\par
It remains to compute the case for $b>1$.
Choose an integer $m$ satisfying the inequality $0<a+bm<b$.
Equivalently, $-\frac{a}{b}<m<1-\frac{a}{b}$.
This can always be achieved, since the fact that $a>0,b>1$ are relatively 
prime integers implies that $b$ does not divide $a$. i.e. 
$\frac{a}{b}$ is not an integer.
Left multiplication by the matrix 
$\begin{pmatrix}
1 & m\\
0 & 1
\end{pmatrix}$
maps the vector $(a,b)$ to $(a+bm,b)$.
Noting that the matrix above is in GL$(2,\mathbb{Z})$, Lemma 
3.1.2 in my PhD thesis \cite{S2} implies that the ring 
$k[\sigma\cap\mathbb{Z}^2]$ 
is isomorphic to the ring $k[\sigma'\cap\mathbb{Z}^2]$, where 
$\sigma'$ is the cone in $\mathbb{R}^2$ 
generated by $e_1,(a+bm)e_1+be_2$, and $0<a+bm<b$ are 
relatively prime integers.
Thus, we may assume that the cone $\sigma$ is generated by 
$e_1,ae_1+be_2$ 
for some positive integers $a<b$, which are relatively prime.
\end{remark}
\begin{Lemma}
Let $\sigma$ be the cone in $\mathbb{R}^2$ generated by 
$e_1,ae_1+be_2$, 
where $a<b$ are relatively prime positive integers.
Let $X$ be the affine toric variety 
$\operatorname{Spec}(k[\sigma\cap\mathbb{Z}^2])$.\\
Then $G_n(X)\cong G_n(k)\oplus\mathbb{Z}/b\mathbb{Z}$ for every 
non-negative even number $n$, and $G_n(X)\cong G_n(k)$ for 
every positive odd number $n$.
\end{Lemma}
\begin{proof}
Note that the points in $\sigma\cap\mathbb{Z}^2$ are all integer-coordinate 
translates of the integer-coordinate points in the region $T$ bounded by 
the non-negative x-axis, the line $y=\frac{b}{a} x$ and the line $x=a$.
Thus, the semigroup generators of $\sigma\cap\mathbb{Z}^2$ can be chosen 
from this region.
For any integer-coordinate point in this region which lies strictly below 
the line $y=\frac{b}{a} x$, any multiple of the point by a non-negative 
integer stays strictly below the line $y=\frac{b}{a} x$.\par
And for any non-negative integer linear combination of a pair of integer-
coordinate points in this region which lie strictly below 
the line $y=\frac{b}{a} x$, the resulting point lies between the two rays 
generated by each of the two points. So it also lies strictly below the 
line $y=\frac{b}{a} x$.
By induction on the number of points used for the non-negative integer 
linear combination, we deduce that the resulting point of any such linear 
combination lies strictly below the line $y=\frac{b}{a} x$.
Therefore, the point $ae_1+be_2$ is not in the semigroup generated by all 
the other integer-coordinate points in the region $T$, since
all those points lie strictly below the line $y=\frac{b}{a} x$.
i.e. $ae_1+be_2$ is in the list of semigroup generators for the 
semigroup $\sigma\cap\mathbb{Z}^2$, which are contained in the 
region $T$.
Let $v=x^ay^b$.
I claim that for any point $(c,d)\in\mathbb{Z}^2$ which lies strictly below 
the line $y=\frac{b}{a} x$ and is contained in the region $T$, the 
quotient image of $x^cy^d$ in $R/xR$ is nilpotent.\par
Since the zero element of the ring $R/xR$ corresponds to the points in 
$\sigma\cap\mathbb{Z}^2$ which lie in the cone $\sigma$ when translated to 
the left by one unit, and for large enough positive integers $l$, the point 
$(cl,dl)$ has this property,the quotient image of $x^{cl}y^{dl}$ in $R/xR$ 
is zero for large enough integers $l>0$.\\ 
So the quotient image of such $x^cy^d$ in $R/xR$ is nilpotent.\par
Therefore, by noting that $(a,b)$ is the only point in the intersection of 
$\sigma\cap\mathbb{Z}^2$ and the region $T$ which is not strictly below the 
line $y=\frac{b}{a} x$, we deduce that the quotient image in $R/xR$ of all 
semigroup generators for $\sigma\cap\mathbb{Z}^2$ (chosen from the region 
$T$), except $ae_1+be_2$, are nilpotent elements of $R/xR$.\\
i.e. $(R/xR)_{red}\cong k[v]\cong k[t]$ is a principal ideal 
domain.\par
Since $R$ is the semigroup algebra $k[\sigma\cap\mathbb{Z}^2]$, the 
elements of the semigroup $\sigma\cap\mathbb{Z}^2$ form a $k$-
basis for the $k$-algebra $R$.\par
Since multiplication by $x$ in the ring $R$ corresponds to right 
translation by one unit for points in $\sigma\cap\mathbb{Z}^2$, there is a  
$k$-basis for $R/xR$ corresponding to the points in $\sigma\cap\mathbb{Z}^2$
which lie outside the cone $\sigma$ when translated one unit to the left.\par
Hence, there is a $k$-basis for the $k$-algebra $R/xR$ consisting 
of elements of the form 
\[
v^l.1,v^l.(xy),v^l.(xy^2),...,v^l.(xy^{\floor{\frac{b}
{a}}}),v^l.(x^2y^{\floor{\frac{b}{a}}+1)}),v^l.(x^2y^{\floor{\frac{b}
{a}}+2}),...,v^l(x^2y^{\floor{2(\frac{b}{a})}}),
\]
...
\[
v^l.(x^{a-1}y^{\floor{(a-2)(\frac{b}{a})}+1}),...,v^l.(x^{a-
1}y^{\floor{(a-1)(\frac{b}{a})}}),v^l(x^ay^{\floor{(a-1)(\frac{b}{a})}+1}),...,v^l.(x^ay^{b-1}), 
\]
where $l=0,1,2,...$.
Therefore, we have 
\begin{multline*}
R/xR=k[v].1\oplus k[v].(xy)...\oplus k[v]
(xy^{\floor{\frac{b}{a}}})
\oplus k[v](x^2y^{\floor{\frac{b}{a}}+1)})...\\\oplus k[v](x^2y^{\floor{2(\frac{b}{a})}})...\oplus k[v](x^{a-1}y^{\floor{(a-2)(\frac{b}{a})}+1})\oplus k[v](x^{a-1}y^{\floor{(a-1)(\frac{b}{a})}+1})\\\oplus k[v](x^ay^{\floor{(a-1)(\frac{b}{a})}+1})...\oplus k[v](x^ay^{b-1}).
\end{multline*}
Therefore, the class $[R/xR]\in G_0(R/xR)$ corresponds to $(1)+(\floor{\frac{b}
{a}})+(\floor{2(\frac{b}{a})}-\floor{\frac{b}{a}})+...+(\floor{(a-1)(\frac{b}{a})})-(\floor{(a-2)
(\frac{b}{a})})+(\floor{a(\frac{b}{a})}-\floor{(a-1)(\frac{b}{a})}-1)=b\in 
\mathbb{Z}\cong G_0((R/xR)_{red})$.\par
By a reasoning similar to the computation of 
$G_0(\operatorname{Spec}k[x,xy,xy^2,...,xy^m])$ as in the proof 
of Lemma 2.2 above, we get that $G_0(X)\cong 
\mathbb{Z}\oplus\mathbb{Z}/b\mathbb{Z}$.
Similar reasoning to the computation of the higher $G$-theory 
of 
$\operatorname{Spec}(k[x,xy,xy^2,...,xy^m])$ yields that 
for every positive even number 
$n$, $G_n(X)\cong G_n(k)\oplus \mathbb{Z}/b\mathbb{Z}$.
And for every positive odd number $n$, $G_n(X)\cong 
G_n(k)$.
\end{proof}
\begin{Theorem}
Let $k$ be an algebraically closed field of characteristic 
zero.\\
Let $\sigma$ be a 2-dimensional strongly convex, rational, 
polyhedral cone in $\mathbb{R}^2$.\\
Let $\delta$ be the determinant of the matrix taking the 
two minimal generators of the cone $\sigma$ as its 
columns.\\
Let $X$ be the affine toric variety $\operatorname{Spec}
(k[\sigma^{\vee}\cap\mathbb{Z}^2])$.\\
Then $G_0(X)\cong\mathbb{Z}\oplus\mathbb{Z}/\delta\mathbb{Z}$.\\
In general, we have $G_n(X)\cong G_n(k)\oplus\mathbb{Z}/\delta\mathbb{Z}$ if 
$n\geq 0$ is even, and $G_n(X)\cong G_n(k)$ if $n\geq 0$ is odd.
\end{Theorem}
\begin{proof}
Recall that there is a matrix $M\in GL(2,\mathbb{Z})$ such 
that left multiplication by $M$ maps the cone generators
of $\sigma$ to $e_1,ae_1+be_2$, where $a,b>0$ are 
relatively prime integers.
And the image of the cone $\sigma$ is the cone $\tau$ in $\mathbb{R}^2$ 
with minimal generators $e_1,ae_1+be_2$.\par
Note that the dual of $\tau$ has minimal generators $e_2,be_1-ae_2$.\\
Since the integers $a,b$ are relatively prime, there exists integers $s,t$
such that $as+bt=1$.\\
Let $A$ be the matrix
$\begin{pmatrix}
b & s\\
-a & t\\
\end{pmatrix}$.
Then $\mathrm{det}(A)=bt+as=1$, so $A\in GL(2,\mathbb{Z})$, and
$A^{-1}$ is the matrix
$\begin{pmatrix}
t & -s\\
a & b\\
\end{pmatrix}$.\par
And left mutiplication by $A^{-1}$ maps $(0,1)$ to $(-s,b)$, maps $(b,-a)$
to $(1,0)$.\\
Since $A^{-1}$ represents an automorphism of the lattice $\mathbb{Z}^2$,
the cone $\tau^{\vee}$ is isomorphic to the cone $\rho$ in $\mathbb{R}^2$ 
generated by $e_1,(-s,b)$.\par
Since $X$ is the affine toric variety $\operatorname{Spec}
(k[\sigma^{\vee}\cap\mathbb{Z}^2])$, by applying Lemma 3.1.2
in \cite{S2}, we obtain that $X$ is isomorphic to 
$\operatorname{Spec}(k[\rho\cap\mathbb{Z}^2])$.\par
Choose an integer $m$ such that $-s+bm>0$.
This is always possible, since we only need $m>\frac{s}{b}$.
Note that since $as+bt=1$, we know that the integers $-s,b$ are 
relatively prime.\par
Let $B$ be the matrix
$\begin{pmatrix}
1 & m\\
0 & 1\\
\end{pmatrix}$ in $GL(2,\mathbb{Z})$.
Left multiplication by the matrix $B$ represents an 
automorphism of the lattice $\mathbb{Z}^2$, 
and this transformation maps $(-s,b)$ to $(-s+bm,b)$ and fixes $(1,0)$.\\
Since $-s,b$ are relatively prime integers, so are $-s+bm$ and $b$.
Hence, by replacing $(-s,b)$ with $(-s+bm,b)$, we may assume that the cone 
generator $(-s,b)$ for $\rho$ satisfies $-s,b>0$ and $-s,b$ are relatively 
prime integers.\par
Choose an integer $n$ such that $-s+bn\in (0,b)$.
This is always possible, since we only need $n\in (\frac{s}{b},\frac{b+s}
{b})$, and that $\frac{s}{b}$ is not an integer.
By replacing $(-s,b)$ with $(-s+bn,b)$, we may assume that $0<-s<b$
are relatively prime.
By Theorem 3.6.6 on the $G$-theory of affine toric surfaces, we 
have
$G_0(X)\cong\mathbb{Z}\oplus\mathbb{Z}/b\mathbb{Z}$.\par
Since the matrix taking the minimal generators of the 
cone $\tau$ as its columns is
$\begin{pmatrix}
1 & a\\
0 & b\\
\end{pmatrix}$, whose determinant is $b$, and this matrix is obtained by 
left multiplication by a finite sequence of matrices in 
$GL(2,\mathbb{Z})$,
we deduce that the determinant $\delta$ of the matrix taking 
the minimal 
generators of the cone $\sigma$ we started with must be $b$ or 
$-b$.\par
Therefore, we conclude that 
$G_0(X)\cong\mathbb{Z}\oplus\mathbb{Z}/\delta\mathbb{Z}$.
Similarly, we conclude that $G_n(X)\cong 
G_n(k)\oplus\mathbb{Z}/\delta\mathbb{Z}$ if 
$n\geq 0$ is even, and $G_n(X)\cong G_n(k)$ if $n\geq 0$ is 
odd.
\end{proof}
Now we prove a lemma, which will be used to compute all the $G$-theory groups of resolution of singularities of the affine 
toric surface $\operatorname{Spec}(k[x,xy,xy^2,...,xy^d])$.\par
\begin{Lemma}
Let $k$ be a field and $X$ be an algebraic variety over $k$.
Let $x\in X$ be a smooth $k$-rational point.
Choose an open subscheme $U$ of $X$ containing $x$.
Then for every non-negative integer $n$, we have 
$G_n(X)\cong G_n(k)\oplus \widetilde{G}_n(X)$ and 
$G_n(U)\cong G_n(k)\oplus \widetilde{G}_n(U)$.
Let $Z=X-U$. Then the $G$-theory localization sequence for $X$ 
and $U$ induces a long exact sequence
\begin{equation*}
...\rightarrow 
G_n(Z)\rightarrow \widetilde{G}_n(X)\rightarrow 
\widetilde{G}_n(U)\rightarrow G_{n-1}(Z)\rightarrow 
...\rightarrow \widetilde{G}_0(U)\rightarrow 0.
\end{equation*}
\end{Lemma}
\begin{proof}
First, we show that that the smooth $k$-rational point $x\in X$
is isomorphic to $\operatorname{Spec}(k)$ as a closed
subscheme of $X$. Since $X$ is a $k$-scheme of finite type, 
$x$ is a closed point of an affine open subset 
$V=\operatorname{Spec}(A)$ of $X$. It corresponds to a 
maximal ideal $m$ of the commutative ring $A$.  Hence, 
$O_{X,x}\cong O_{V,x}\cong A_m$ and $m_x\cong mA_m$.
Therefore, the residue field $\kappa(x)=O_{X,x}/m_x\cong 
A_m/mA_m\cong A/m$.
Since $x$ is a $k$-rational point of $X$, $\kappa(x)\cong k$.
Therefore, $x=V(m)\cong \operatorname{Spec}(A/m)\cong 
\operatorname{Spec}\kappa(x)\cong \operatorname{Spec}(k)$.
Next, we construct the direct sum decompositions 
$G_n(X)\cong G_n(k)\oplus \widetilde{G}_n(X)$ and 
$G_n(U)\cong G_n(k)\oplus \widetilde{G}_n(U)$ for every non-negative integer $n$.
Note that the identity morphism on the $k$-rational point 
$x\in X$ is the composition of the inclusion of $x$ into the 
smooth locus of $X$, denoted by $X_{smooth}$, followed by 
the inclusion of $X_{smooth}$ into $X$ and then by the 
natural projection $\pi:X\rightarrow \operatorname{Spec}(k)$.
Since a morphism into a regular scheme of finite Krull 
dimension is always a morphism of finite Tor-dimension, the 
inclusion of $x$ into $X_{smooth}$ and the projection 
$\pi:X\rightarrow \operatorname{Spec}(k)$ are both morphisms 
of finite Tor-dimension.
Since the inclusion of $X_{smooth}$ into $X$ is an open 
immersion, it is flat and hence also a morphism of finite 
Tor-dimension.
Since the $G$-theory of Noetherian schemes is contravariant 
with respect to morphisms of finite Tor-dimension, we have 
$f^*\circ \pi^*=id$. Here $f:\operatorname{Spec}
(k)\rightarrow X$ is the inclusion of $x$ into $X$, which 
decomposes into the inclusion of $x$ into $X_{smooth}$, 
followed by the inclusion of $X_{smooth}$ into $X$.
Therefore, the homomorphisms $\pi^*:G_n(k)\rightarrow 
G_n(X)$ induced by the projection $\pi$, are split 
injections for every non-negative integer $n$. 
Define $\widetilde{G}_n(X)$ as the cokernel of 
$\pi^*:G_n(k)\rightarrow G_n(X)$ for every non-negative 
integer $n$. 
Then we have a split short exact sequence for every non-negative integer $n$:
\begin{equation*}
0\rightarrow G_n(k)\rightarrow G_n(X)\rightarrow 
\widetilde{G}_n(X)\rightarrow 0
\end{equation*}
We define $\widetilde{G}_n(U)$ similarly.
Now we construct the induced long exact sequence 
\begin{equation*}
...\rightarrow 
G_n(Z)\rightarrow \widetilde{G}_n(X)\rightarrow 
\widetilde{G}_n(U)\rightarrow G_{n-1}(Z)\rightarrow 
...\rightarrow \widetilde{G}_0(U)\rightarrow 0
\end{equation*}
Define $G_n(Z)\rightarrow \widetilde{G}_n(X)$ as the composition
$G_n(Z)\rightarrow G_n(X)\rightarrow \widetilde{G}_n(X)$, where 
the first homomorphism is part of the $G$-theory 
localization sequence for $X$, $U$ and $Z=X-U$.
Define $\widetilde{G}_n(X)\rightarrow \widetilde{G}_n(U)$ by 
the 
universal property of cokernels in the abelian category of 
abelian groups. They are induced by the homomorphisms 
$G_n(X)\rightarrow G_n(U)$, which are part of the $G$-theory 
localization sequence for $X$, $U$ and $Z=X-U$.
Define $\widetilde{G}_n(U)\rightarrow G_{n-1}(Z)$ by pulling 
back elements along the surjection $G_n(U)\rightarrow 
\widetilde{G}_n(U)$ and then pushing them along 
$G_n(U)\rightarrow G_{n-1}(Z)$.
A diagram chase verifies that the homomorphisms 
$\widetilde{G}_n(U)\rightarrow G_{n-1}(Z)$ are well-defined and 
fit into a long exact sequence 
\begin{equation*}
...\rightarrow 
G_n(Z)\rightarrow \widetilde{G}_n(X)\rightarrow 
\widetilde{G}_n(U)\rightarrow G_{n-1}(Z)\rightarrow 
...\rightarrow \widetilde{G}_0(U)\rightarrow 0
\end{equation*}
as required.
\end{proof}
\begin{Theorem}
Let $k$ be a field and let $d>0$ be an integer. Let $\sigma$ be 
the cone in $\mathbb{R}^2$ generated by $e_2,de_1-e_2$.
Let $X=\operatorname{Spec}(k[\sigma^{\vee}\cap\mathbb{Z}^2])$ 
be 
the affine toric surface over the field $k$ associated to the 
cone $\sigma$. Let $\tilde{X}$ be the resolution of 
singularities of $X$ obtained by adding the ray in 
$\mathbb{R}^2$ generated by $e_1$ to the fan of $X$.
Then $G_n(\tilde{X})\cong G_n(k)^{\oplus 2}$ for all $n\geq 0$.
\end{Theorem}
\begin{proof}
Let $\sigma_1$ be the cone in $\mathbb{R}^2$ generated by 
$e_2, e_1$, $\sigma_2$ be the cone in $\mathbb{R}^2$ generated 
by $e_1, de_1-e_2$. And let $\tau$ be the cone in $\mathbb{R}^2$ 
generated by $e_1$. Let $U_{\sigma_i}=\operatorname{Spec}(k[\sigma_i^{\vee}\cap\mathbb{Z}^2])$ be the affine toric variety over the field $k$ associated to the cones $\sigma_i$, $i=1, 2$. And let $U_{\tau}=\operatorname{Spec}(k[\tau^{\vee}\cap\mathbb{Z}^2])$ be the affine toric variety 
over the field $k$ associated to the cone $\tau$.
Then $U_{\tau}=U_{\sigma_1}\cap\ U_{\sigma_2}$,
$G_n(U_{\sigma_1})\cong G_n(k)$ and $G_n(U_{\sigma_2})\cong G_n(k)$ for all $n\geq 0$.
$U_{\tau}\cong\operatorname{Spec}(k[x,y,y^{-1}])$.
By Fundamental Theorem of $G$-theory of Noetherian rings, $G_n(U_{\tau})\cong G_n(k)\oplus G_{n-1}(k)$ for all $n\geq 0$.\\
Since $\{U_{\sigma_1}, U_{\sigma_2}\}$ forms an affine open 
cover of $\tilde{X}$, there is a homotopy Cartesian square of 
$G$-theory spectra:
\begin{tikzcd}
G(\tilde{X}) \arrow{r} \arrow{d}
& G(U_{\sigma_1}) \arrow {d}\\
G(U_{\sigma_2}) \arrow {r}
& G(U_{\tau})
\end{tikzcd}.\\
Note that for every field extension $k'/k$, the corresponding 
base change $U_{\tau,k'}:=U_{\tau}\times_k \operatorname{Spec}(k')$ is a regular $k'$-scheme.
Hence, the $k$-scheme $U_{\tau}$ is geometrically regular.
$U_{\tau}$ is clearly a $k$-scheme of finite type, so it is a 
smooth $k$-scheme. 
The set of $k$-rational points of $U_{\tau}$ is $U_{\tau}(k)=\mathrm{Hom}_{k-alg}(k[x,y,y^{-1}],k)\cong k\times k^{\times}$.
Hence, the set of smooth $k$-rational points of $U_{\tau}$ is 
$k\times k^{\times}$, which is not empty.
By Lemma 2.7 above, if we set $Z=X-U_{\sigma_1}$, 
then the homotopy Cartesian square of $G$-theory spectra above 
induces a long exact sequence 
\begin{align*}
...\rightarrow \tilde{G}_n(\tilde{X})\rightarrow \tilde{G}_n(U_{\sigma_1})\oplus \tilde{G}_n(U_{\sigma_2})\rightarrow \tilde{G}_n(U_{\tau})\rightarrow \tilde{G}_{n-1}(\tilde{X})\rightarrow ...
\end{align*}
Since $\tilde{G}_n(U_{\sigma_1})\cong G_n(U_{\sigma_1})/G_n(k)\cong 0$ and $\tilde{G}_n(U_{\sigma_1})\cong G_n(U_{\sigma_2})/G_n(k)\cong 0$
for all $n\geq 0$, the homomorphism $\tilde{G}_n(U_{\tau})\rightarrow \tilde{G}_{n-1}(\tilde{X})$ is 
always an isomorphism.
Hence, $G_n(\tilde{X})\cong\tilde{G}_n(\tilde{X})\oplus G_n(k)\cong\tilde{G}_{n+1}(U_{\tau})\oplus G_n(k)\cong [G_{n+1}(U_{\tau})/G_{n+1}(k)]\oplus G_n(k)\cong G_n(k)^{\oplus 2}$
for all $n\geq 0$.
\end{proof}

Now we prove a special case of the conjecture in \cite{S1}.\par
\begin{Theorem}
Let $X$ be any affine, smooth toric variety over a field $k$.
Then the Chow group of codimension 2 cycles $A^2(X)$ is a finite 
abelian group whose order divides $|\delta|$, where $\delta$ is 
defined above. In fact, the Chow group $A^2(X)$ is trivial here.
\end{Theorem}
\begin{proof}
Since $X$ is an affine, smooth toric variety over the field $k$, 
it has the form $\operatorname{Spec}(k[\sigma^{\vee}\cap\mathbb{Z}^n])$
for some smooth $r$-dimensional cone $\sigma$ in $\mathbb{R}^n$
($0\leq r\leq n$).
Here $n$ is the Krull dimension of $X$.
Let $u_1,u_2,...,u_r$ be the minimal generators of the smooth 
cone $\sigma$.
Since the minimal generators $u_1,u_2,...,u_r$ of the smooth 
cone $\sigma$ extend to a $\mathbb{Z}$-basis 
$\{u_1,u_2,...,u_r,u_{r+1},...,u_n\}$ for the lattice 
$\mathbb{Z}^n$, there exists a matrix $A\in GL(n,\mathbb{Z})$ 
such that $Au_i=e_i$ for $i=1,2,...,n$, where $\{e_1,e_2,...,e_n\}$ is
the standard $\mathbb{Z}$-basis for $\mathbb{Z}^n$.
By Lemma 3.2 in \cite{S1}, there is a ring isomorphism 
$k[\sigma^{\vee}\cap\mathbb{Z}^n]\rightarrow k[\tau^{\vee}\cap\mathbb{Z}^n]$, where $\tau$ is the $r$-dimensional smooth cone in $\mathbb{R}^n$ generated by $e_1,e_2,...,e_r$.
Since $k[\tau^{\vee}\cap\mathbb{Z}^n]\cong k[t_1,t_2,...,t_r,t_{r+1}^{\pm},...,t_n^{\pm}]$, 
where $k[t_1,t_2,...,t_r,t_{r+1}^{\pm},...,t_n^{\pm}]$ denotes 
the localization of the polynomial ring in $n$ indeterminates 
$t_1,t_2,...,t_n$ over the field $k$ at the element 
$t_{r+1}t_{r+2}...t_n$, the ring 
$k[\sigma^{\vee}\cap\mathbb{Z}^n]$ is isomorphic to this 
localization. 
i.e., $X=\operatorname{Spec}(k[\sigma^{\vee}\cap\mathbb{Z}^n])\cong\operatorname{Spec}(k[t_1,t_2,...,t_r,t_{r+1}^{\pm},...,t_n^{\pm}])\cong \mathbb{A}^r_k\times\mathbb{G}_{m,k}^{n-r}$.
Therefore, $A^2(X)=A_{n-2}(X)\cong A_{n-2-r}(\mathbb{G}_{m,k}^{n-r})\cong A^2(\mathbb{G}_{m,k}^{n-r})$, by the homotopy property of Chow groups.
If $r=n$, then $X=\mathbb{A}^n_k$, so $A^2(X)=0$ as required.
Now assume $r<n$. Let $N=n-r$. We show that $A^2(\mathbb{G}_{m,k}^N)=0$ by induction on $N$.
If $N=1$, then $A^2(X)\cong A^2(\mathbb{G}_{m,k})=0$ since
$\mathbb{G}_{m,k}$ is a 1-dimensional variety.
Suppose that for some integer $N\geq 2$, $A^2(\mathbb{G}_{m,k}^i)=0$ has been proved for all positive 
integers $i<N$.
By the Künneth formula for Chow groups of linear schemes in 
\cite{T}, if $Y, Z$ are two toric varieties over the field $k$, then there is an isomorphism $A_{*}(Y\times Z)\cong A_{*}(Y)\otimes_{\mathbb{Z}} A_{*}(Z)$.\\
Hence, $A^2(\mathbb{G}_{m,k}^N)$ is isomorphic to the direct sum
\begin{align*}
[A_0(\mathbb{G}_{m,k})\otimes_{\mathbb{Z}}A_{N-2}(\mathbb{G}_{m,k}^{N-1})]\oplus [A_1(\mathbb{G}_{m,k})\otimes_{\mathbb{Z}}A_{N-3}(\mathbb{G}_{m,k}^{N-1})]\oplus...\oplus[A_{N-2}(\mathbb{G}_{m,k})\otimes_{\mathbb{Z}}A_0(\mathbb{G}_{m,k}^{N-1})].
\end{align*}
Since $k[x,x^{-1}]$ is a Noetherian UFD,
$A_0(\mathbb{G}_{m,k})\cong\mathrm{Cl}(\mathbb{G}_{m,k})=0$.
Note that $A_{N-3}(\mathbb{G}_{m,k}^{N-1})\cong A^2(\mathbb{G}_{m,k}^{N-1})$, which vanishes by the induction 
hypothesis.
And for all integers $i>1$, $A_i(\mathbb{G}_{m,k})=0$ since 
$\mathbb{G}_{m,k}$ is a 1-dimensional variety.
Therefore, each direct summand above vanishes, completing
the induction on $N$.
\end{proof}
\section{$G_0$ with rational coefficients of complete, 
simplicial toric varieties}
\begin{Theorem}
Let $X$ be a complete, simplicial toric variety over an 
algebraically closed field of characteristic zero. Then the 
dimension of the rational vector space $G_0(X)\otimes\mathbb{Q}$
is the sum of the Betti numbers of even degrees of $X$. 
The Betti numbers here are defined to be $b_{2k}(X_{\Sigma})=\mathrm{dim}H^{2k}(X_{\Sigma},\mathbb{Q})$.
This can be expressed in terms of the number of cones of 
various dimensions in the fan of $X$.
\end{Theorem}
\begin{proof}
By the Lefschetz principle, it suffices to assume that the base 
field of $X$ is $\mathbb{C}$.
By Fulton's Riemann-Roch theorem for algebraic schemes in 
\cite{F}, we have
$G_0(X)\otimes \mathbb{Q}\cong A_*(X)\otimes \mathbb{Q}$.\\
Since the fan $\Sigma$ of $X$ is complete and simplicial, we 
have $A_*(X)\otimes \mathbb{Q}\cong A^*
(X)\otimes\mathbb{Q}$.\\
And by Theorem 12.5.3 in \cite{CLS}, we have $A^*(X)\otimes 
\mathbb{Q}\cong H^*(X,\mathbb{Q})$, where $H^*(X,\mathbb{Q})$ 
denotes the 
singular cohomology ring of $X$ with rational coefficients.\\
By Theorem 12.3.11 in \cite{CLS} above, the singular cohomology 
groups of $X$ with rational coefficients vanish in odd degrees.\\
Hence, the dimension of the rational vector 
space $G_0(X)\otimes \mathbb{Q}$ is equal to the sum of the 
Betti numbers of $X$ of even degree.\\
By Theorem 12.3.12 in \cite{CLS}, the Betti numbers of $X$ of 
even degrees are given by 
\[
b_{2k}=\sum_{i=k}^{n} (-1)^{i-k}\binom{i}{k} |\Sigma(n-
i)|
\] 
for all $k$.\\
And these Betti numbers satisfy $b_{2k}=b_{2n-2k}$.
\end{proof}
\section{Higher $G$-theory of weighted projective spaces}
Now we compute all the $G$-theory groups of any weighted 
projective space over an arbitrary field. We prove a lemma 
first, which allows us to apply Lemma 2.7 above to give our 
main result in this section.\par
\begin{Lemma}
Let $k$ be a field and let $X$ be a toric variety over $k$ 
associated to a fan $\Sigma$ in $\mathbb{R}^n$.
Then the set of all smooth $k$-rational points of $X$ is non-empty. 
\end{Lemma}
\begin{proof}
Since the smooth locus of a geometrically reduced $k$-variety 
is open and dense in the variety \cite{P} and the 
$k$-variety $X$ is geometrically reduced, its smooth locus 
is non-empty.
By Theorem 2.6.20 above, the smooth locus of the toric variety 
$X$ is the union of 
the affine toric varieties $U_{\sigma}:=\operatorname{Spec}
(k[\sigma^{\vee}\cap\mathbb{Z}^n])$ associated to the 
smooth cones $\sigma$ in the fan of $X$. Here 
$n=\mathrm{dim}(X)$. 
Let $\sigma$ be any $r$-dimensional smooth cone in the fan 
of $X$. Its set of minimal generators $\{u_1,u_2,...,u_r\}$
can be extended to a $\mathbb{Z}$-basis $\{u_1,u_2,...,u_r,u_{r+1},...,u_n\}$ for $\mathbb{Z}^n$.
Let $\phi:\mathbb{Z}^n\rightarrow \mathbb{Z}^n$ be defined 
by setting $\phi(e_i)=u_i$ for $i=1,2,...,n$ and extended 
by $\mathbb{Z}$-linearity to $\mathbb{Z}^n$, where $\{e_1,e_2,..., e_n\}$ is the standard $\mathbb{Z}$-basis for 
$\mathbb{Z}^n$.
Then $\phi$ is a transformation in $GL(n,\mathbb{Z})$, 
since it maps a $\mathbb{Z}$-basis of $\mathbb{Z}^n$ to 
another one.
By Lemma 3.1.3 above, if $A$ denotes the matrix in 
$GL(n,\mathbb{Z})$ such that $Ae_i=u_i$ for $i=1,2,...,n$, 
then left multiplication by $(A^{-1})^t$ induces a ring 
isomorphism 
$k[x_1,x_2,...,x_r,x_{r+1}^{\pm},...,x_n^{\pm}]\rightarrow 
k[\sigma^{\vee}\cap\mathbb{Z}^n]$.
Hence, we have $U_{\sigma}\cong\operatorname{Spec}(k[x_1,x_2,...,x_r,x_{r+1}^{\pm},...,x_n^{\pm}])$, so that
$U_{\sigma}\cong\mathbb{A}_k^r\times_k \mathbb{G}_m^{n-r}$.
By Proposition 2.6.21 above, if we work in the category of 
schemes of 
finite type over the field $k$ and set 
$S=\operatorname{Spec}(k)$, then the set of 
$k$-rational points of the fiber product 
$\mathbb{A}_k^r\times_k \mathbb{G}_m^{n-r}$ is the 
Cartesian product of the $k$-rational points of the factors.
i.e., $U_{\sigma}(k)=\mathbb{A}^r_k(k)\times 
\mathbb{G}_m^{n-r}(k)=k^r\times (k^*)^{n-r}$.
Since $U_{\sigma}$ is an affine open subset of the smooth 
locus of toric variety $X$, the set of $k$-rational 
points of the smooth locus of $X$ contains the set of $k$-rational points of
$U_{\sigma}$, so the set of smooth $k$-rational points 
of $X$ is non-empty.
\end{proof}

\begin{Theorem}
Let $X$ be any $d$-dimensional weighted projective space over 
an arbitrary field $k$. i.e. $X=\operatorname{Proj}(S)$ for 
some graded ring $S=k[x_0,x_1,...,x_d]$, where $x_0,x_1,...,x_d$
are indeterminates of degrees $a_0,a_1,...,a_d$ for some 
relatively prime positive integers $a_0,a_1,...,a_d$.
Then $G_n(X)\cong G_n(k)^{\oplus 
(d+1)}$ for all $n$.
\end{Theorem}
\begin{proof}
We proceed by induction on the dimension $d$ of $X$. When $d=1$, 
$X$ is the weighted projective line over the field $k$.
It is isomorphic to the projective line over $k$.
By projective bundle formula, $G_n(X)\cong G_n(k)^{\oplus 2}$ 
for every non-negative integer $n$. Hence, the base case is 
proved.\par
Now consider the case $d=2$, so that $X$ has the form $\mathbb{P}
(a_0,a_1,a_2)$ for some relatively prime positive integers
$a_0,a_1,a_2$. Set $S$ to be the polynomial ring $k[x_0,x_1,x_2]$
where each $x_i$ has degree $a_i$. Then $X=\operatorname{Proj}
(S)$. Set $Z=V_{+}(x_2)$ and $U=D_{+}(x_2)$.
Then there is a long exact sequence induced by the $G$-theory 
localization sequence for the Noetherian scheme $X$. 
It has the form
\[
...\rightarrow G_n(Z)\rightarrow 
\widetilde{G}_n(X)\rightarrow \widetilde{G}_n(U)\rightarrow G_{n-1}
(Z)\rightarrow ...\rightarrow 0 (*)
\]
Here $G_n(X)\cong \widetilde{G}_n(X)\oplus G_n(k)$ for every non-
negative integer $n$, same for $U$, see Lemma 2.7 above.
Note that $U$ is isomorphic to the scheme $\operatorname{Spec}(S_{(x_2)})$.
Every element of the graded localization $S_{(x_2)}$ has the 
form $\frac{f}{x_2^n}$, where $f\in S$ is homogeneous and has 
the same degree as the monomial $x_2^n$. Since $x_2$ has degree 
$a_2$, every monomial 
$x_0^{\alpha_0}x_1^{\alpha_1}x_2^{\alpha_2}$ appearing in the 
polynomial $f$ has degree $na_2$,
i.e., $na_2=a_0\alpha_0+a_1\alpha_1+a_2\alpha_2$.
Hence, the quotient 
$\frac{x_0^{\alpha_0}x_1^{\alpha_1}x_2^{\alpha_2}}{x_2^n}$ 
simplifies into 
\[(\frac{x_0}{{x_2}^{\frac{a_0}
{a_2}}})^{\alpha_0}(\frac{x_1}{{x_2}^{\frac{a_1}
{a_2}}})^{\alpha_1}.
\]
Since $\frac{x_0}{{x_2}^{\frac{a_0}{a_2}}}$ and $\frac{x_1}
{{x_2}^{\frac{a_1}{a_2}}}$ are algebraically independent over 
the field $k$, the ring $S_{(x_2)}$ is isomorphic to the 
polynomial ring in two variables over the field $k$.
Therefore, $U$ is isomorphic to the scheme $\mathbb{A}_k^2$.
Thus, $\widetilde{G}_n(U)=0$ for all non-negative integers $n$.
By the long exact sequence $(*)$ above, the homomorphism 
$G_n(Z)\rightarrow \widetilde{G}_n(X)$ is an isomorphism for 
every 
non-negative integer $n$. 
Note that $Z$ is the weighted projective line $\mathbb{P}
(a_0,a_1)$. By the base step for induction above, we know that
$G_n(Z)\cong G_n(k)^{\oplus 2}$ for every non-negative integer 
$n$. Therefore, we have $G_n(X)\cong \widetilde{G}_n(X)\oplus 
G_n(k)\cong G_n(k)^{\oplus 3}$ for every non-negative integer 
$n$.\par
Now assume that $d>2$. Then $X$ has the form 
$\operatorname{Proj}(S)$, where $S$ is the polynomial ring 
$k[x_0,x_1,...,x_d]$ such 
that each variable $x_i$ has degree $a_i$, and 
$a_0,a_1,...,a_d$ 
are relatively prime positive integers. Similar to the case 
$d=2$, we can show that the ring $S_{(x_d)}$ is the polynomial 
ring over the field $k$ with $d$ variables.
By setting $Z=V_{+}(x_d)$ and $U=D_{+}(x_d)$ and applying the 
long exact sequence $(*)$ induced by the $G$-theory 
localization 
sequence of $X$, we conclude that $G_n(Z)\rightarrow 
\widetilde{G}_n(X)$ is always an isomorphism.
Note that $Z$ is isomorphic to the $(d-1)$-dimensional weighted 
projective space $\mathbb{P}(a_0,a_1,...,a_{d-1})$.
By the induction hypothesis, $G_n(Z)\cong G_n(k)^{\oplus d}$ 
for every non-negative integer $n$.
Hence, we conclude that $G_n(X)\cong \widetilde{G}_n(X)\oplus 
G_n(k)\cong G_n(k)^{\oplus (d+1)}$ for every non-negative 
integer $n$.
The proof is complete.
\end{proof}

\section{$G_0,G_1,G_2$ of the product of two weighted 
projective spaces over a field}
\begin{Theorem}
Let $k$ be any field and let $X,Y$ be any two weighted 
projective spaces over the field $k$. Then for $n=0,1,2$,
we have $G_n(X\times_k Y)\cong\bigoplus_{i+j=n}[G_{i}(X)\otimes _{\mathbb{Z}} G_j(Y)]$.
\end{Theorem}
\begin{proof}
Let $S=\operatorname{Spec}(k)$. Then $S$ is a smooth scheme over 
the field $k$, and $X,Y$ are both projective, flat $S$-schemes. 
Since every weighted projective space over a field is a toric 
variety and every toric variety over the field $k$ is a linear 
$k$-scheme \cite{J}, the Künneth spectral sequence for higher 
algebraic K-theory from Theorem 4.2 of \cite{J} applies to the 
pair of $S$-schemes $X,Y$. Let $E^*_{*,*}$ denote this spectral 
sequence. \par
The case $n=0$ is just Corollary 4.3 of \cite{J}. \par
Consider the case $n=1$. We show that $E^{\infty}_{p,0}=0$ for 
all $p>0$ first. 
Recall that the Künneth spectral sequence for 
higher algebraic K-theory is 
\begin{align*}
E^2_{s,t}=\mathrm{Tor}_{s,t}^{\pi_*K(S)}(\pi_*G(X),\pi_*G(Y))\Rightarrow \pi_{s+t}G(X\times_S Y)
\end{align*}
Here $K(S)$ denotes the $K$-theory spectrum of the scheme $S$ 
and $G(X),G(Y)$ denote the $G$-theory spectra of the Noetherian 
schemes $X,Y$ respectively. 
Since $K_0(S)=K_0(k)\cong\mathbb{Z}$, we have 
$E^2_{p,0}=\mathrm{Tor}_p^\mathbb{Z}(G_0(X),G_0(Y))=0$ for all 
$p>0$, since $G_0(X)\cong G_0(k)^{\oplus (1+\mathrm{dim}(X))}$
is a free abelian group. Since $E^{\infty}_{p,0}$ is always a 
subquotient of $E^2_{p,0}$, we have $E^{\infty}_{p,0}=0$ for all 
$p>0$. Since the Künneth spectral sequence for 
higher algebraic K-theory is a convergent first-quadrant 
homological spectral sequence, there is an increasing, finite 
filtration on its abutment. 
Let $F_iG_n(X\times_S Y)$ denote the $i$-th piece of the 
filtration on the abutment.
Then $E^{\infty}_{p,0}=0$ for all $p>0$ implies that $F_pG_p(X\times_S Y)=F_{p-1}G_p(X\times_S Y)$ for all $p>0$.
In particular, we have $F_1G_1(X\times_S Y)=F_0G_1(X\times_S Y)$.
Since $E^2_{s,t}=0$ if either $s<0$ or $t<0$, $E^2_{s,t}=0$ if 
$s+t=1$ and $s>0$. Therefore, $E^{\infty}_{s,t}=0$ if $s+t=1$ 
and $s>0$. Hence, we have $F_0G_1(X\times_S Y)=F_1G_1(X\times_S Y)=F_2G_1(X\times_S Y)=...=G_1(X\times_S Y)$.
Since $E^2_{s,t}=0$ if either $s<0$ or $t<0$, 
$E^{\infty}_{s,t}=0$ if either $s<0$ or $t<0$.
In particular,$F_pG_1(X\times_S Y)/F_{p-1}G_1(X\times_S Y))\cong E^{\infty}_{p,1-p}=0$ for all $p<0$.
i.e., $F_{-1}G_1(X\times_S Y)=F_{-2}G_1(X\times_S Y)=...=0$.
Therefore, $G_1(X\times_S Y)\cong F_0G_1(X\times_S Y)/F_{-1}G_1(X\times_S Y)\cong E^{\infty}_{0,1}$.
Since $E^2_{-2,2}=0=E^2_{2,0}$, $E^2_{0,1}\cong E^3_{0,1}$.
In general, for every integer $r\geq 2$, we have $E^r_{-r,r}=0$
since $E^2_{-r,r}=0$ and $E^r_{-r,r}$ is a subquotient of $E^2_{-r,r}$. Similarly, for every integer $r\geq 2$, $E^r_{r,2-r}=0$.
Hence, the differentials mapping into and out of $E^r_{0,1}$ 
both vanish, so $E^r_{0,1}\cong E^{r+1}_{0,1}$ for all integers
$r\geq 2$. i.e., $E^{\infty}_{0,1}\cong E^2_{0,1}$.
Therefore, $G_1(X\times_S Y)\cong E^{\infty}_{0,1}\cong E^2_{0,1}$, so that 
\begin{align*}
G_1(X\times_S Y)\cong\mathrm{Tor}_{0,1}^{\pi_* K(S)}(\pi_* G(X),\pi_* G(Y))\cong [G_0(X)\otimes_{\mathbb{Z}} G_1(Y)]\oplus [G_1(X)\otimes_{\mathbb{Z}} G_0(Y)].
\end{align*}
Now we compute $G_2(X\times_S Y)$.\\
We show that $F_{p+1}G_2(X\times_S Y)/F_pG_2(X\times_S Y)=0$ for 
all $p\geq 0$.
Since $E^2_{1,1}=\mathrm{Tor}_{1,1}^{\pi_* K(S)}(\pi_* G(X),\pi_* G(Y))=\mathrm{Tor}_1^{\mathbb{Z}}(G_0(X),G_1(Y))\oplus\mathrm{Tor}_1^{\mathbb{Z}}(G_1(X),G_0(Y))$,
and $G_0(X),G_0(Y)$ are both free abelian groups, $E^2_{1,1}=0$.
Since $E^{\infty}_{1,1}$ is a subquotient of $E^2_{1,1}$,
$E^{\infty}_{1,1}=0$. Hence, $F_1G_2(X\times_S Y)/F_0G_2(X\times_S Y)\cong E^{\infty}_{1,1}=0$.
Recall that we showed $E^{\infty}_{p,0}=0$ for all $p>0$.
In particular, $F_2G_2(X\times_S Y)/F_1G_2(X\times_S Y)\cong E^{\infty}_{2,0}=0$.
Since $E^{\infty}_{s,t}=0$ for all $t<0$, we have $F_{p+1}G_2(X\times_S Y)/F_pG_2(X\times_S Y)\cong E^{\infty}_{p+1,2-p-1}=0$ for all $p>1$.\\
i.e., $F_0G_2(X\times_S Y)=F_1G_2(X\times_S Y)=F_2G_2(X\times_S Y)=...=G_2(X\times_S Y)$.
Since $E^{\infty}_{p,2-p}=0$ for all $p<0$, $F_pG_2(X\times_S Y)/F_{p-1}G_2(X\times_S Y)\cong E^{\infty}_{p,2-p}=0$ for all $p<0$.
i.e., $F_{-1}G_2(X\times_S Y)=F_{-2}G_2(X\times_S Y)=F_{-3}G_2(X\times_S Y)=...=0$.
Therefore, $G_2(X\times_S Y)\cong F_0G_2(X\times_S Y)/F_{-1}G_2(X\times_S Y)\cong E^{\infty}_{0,2}$.\\
For every $r\geq 2$, $E^2_{-r,r+1}=0=E^2_{r,2-r+1}$.
Note that 
\begin{align*}
E^2_{2,1}=\mathrm{Tor}_{2,1}^{\pi_* K(S)}(\pi_* G(X),\pi_* G(Y))=\mathrm{Tor}_2^{\mathbb{Z}}(G_0(X),G_1(Y))\oplus\mathrm{Tor}^{\mathbb{Z}}_2(G_1(X),G_0(Y))=0\end{align*}, since $G_0(X),G_0(Y)$ are free abelian groups.\\
Hence, for every $r\geq 2$, $E^r_{-r,r+1}=0=E^r_{r,2-r+1}$.\\
The differentials entering and leaving $E^r_{0,2}$ both vanish 
for all $r\geq 2$, so that $E^2_{0,2}\cong E^{\infty}_{0,2}$.
Therefore, 
\begin{align*}
G_2(X\times_S Y)\cong [G_0(X)\otimes_{\mathbb{Z}} G_2(Y)]\oplus [G_1(X)\otimes_{\mathbb{Z}} G_1(Y)]\oplus [G_2(X)\otimes_{\mathbb{Z}} G_0(Y)].
\end{align*}
\end{proof}

\end{document}